\documentclass[12pt,a4paper]{article}
\usepackage{amsmath}

\numberwithin{equation}{section}
\usepackage{amsthm,amsfonts,amssymb,layout,indentfirst}
\usepackage{float}
\usepackage{multicol}
\usepackage{multirow}
\usepackage{graphicx}
\usepackage{subcaption}
\usepackage{epstopdf}
\usepackage{epsfig}
\usepackage{placeins}
\usepackage{xcolor}
\newlength{\defbaselineskip}
\newcommand{\setlinespacing}[1]%
{\setlength{\baselineskip}{#1 \defbaselineskip}}

\theoremstyle{plain}

\newtheorem*{corollary*}{Corollary}
\newtheorem*{acknowledgement*}{Acknowledgement}

\newtheorem*{remark*}{Remark}

\newtheoremstyle{label}
{3pt}
{3pt}
{\upshape}
{}
{\bfseries}
{}
{.5em}
{}
\theoremstyle{label}
\newtheorem*{lemma*}{Lemma}
\newtheorem*{theorem*}{Theorem}
\newtheoremstyle{citing}
{3pt}
{3pt}
{\upshape}
{}
{\bfseries}
{}
{.5em}
{\thmnote{#3}}
\theoremstyle{citing}

\newtheoremstyle{break}
{9pt}
{9pt}
{\upshape}
{}
{\bfseries}
{.}
{\newline}
{}
\theoremstyle{break}

\title{Qualitative analysis, chaotic structure and exact solution of the  nonlinear seventh-order Caudrey-Dodd-Gibbon-KP equation}

\begin{document}
	\date{}
	\maketitle
	\begin{center}
		\author{\textbf{S. S. Samanta\textsuperscript{1}, S. Sahoo\textsuperscript{2,*}, Vijil Kumar\textsuperscript{3}, Rajib Mia\textsuperscript{4}}
			\\{\textsuperscript{1,2,3} School of Applied Sciences,\\ Kalinga Institute of Industrial Technology,\\ Deemed to be University,\\ Bhubaneswar, Odisha-751024, India,\\
			Email: {\textsuperscript{1}singhsamantasovana@gmail.com}\\
	{\textsuperscript{2,*}subha.bapi25@gmail.com}\\
{\textsuperscript{3}vijilchoudhary@gmail.com}\\
{\textsuperscript{4}rajibmia.90@gmail.com}	}}
		
	\end{center}
	\begin{abstract}
		The main objective of this work is to investigate the traveling wave solution and dynamic characteristics of the (2 + 1)-dimensional seventh-order Caudrey-Dodd-Gibbon-KP (sCDG-KP) equation. Applying the ($\frac{G^{\prime}}{G^{\prime}+G+A}$) method, we examine the exact solution of the (2 + 1)-dimensional seventh-order Caudrey-Dodd-Gibbon-KP (sCDG-KP) equation by altering it into a reduced ODE via a suitable wave transformation. Graphical representations, such as 2D, 3D, and a heat map of the ascertained solution, are present to facilitate comprehension of the empirical relevance of the obtained solutions. As a result, we acquired a bright and anti-kink soliton solution. Next, we alter the ODE into a 2D system of equations to analyze the dynamical behavior of the reduced system via bifurcation analysis, phase portrait, and attractor analysis. During this process, we portray the graphical visualization of the bifurcation phase portrait, 2D phase portrait, 3D phase portrait, time series, chaotic attractor, sensitive analysis, fractal dimension, recurrence plot, and power spectrum of the dynamical system.
		
	\end{abstract}
	\section*{Keywords}
	(2 + 1)-dimensional seventh-order Caudrey-Dodd-Gibbon-KP (sCDG-KP) equation , ($\frac{G^{\prime}}{G^{\prime}+G+A}$) method, Dynamical analysis, Attractor analysis. 
	\section{Introduction}\
	Non-linear partial differential equations have a broad spectrum of uses in many fields such as physics, biology, engineering, economics, non-linear optics and chemical kinematics, etc. Obtaining the exact solution and analyzing the dynamical characteristics of the NLPDEs is becoming very interesting among many researchers due to its wide range of applications to demonstrate the behavior of non-linear complex phenomena. \\
	
	A variety of methods have been established for obtaining the exact solution of an NPDE because these result revels the new features of the wave performance, such methods are: Homotopy perturbation method\cite{Liao2005,Hemeda2012}, Hirota Direct method \cite{Wazwaz2007}, Homotopy analysis method\cite{Matinfar2014,Abbasbandy2008}, BB Acklund transformations method\cite{Khater2002}, Lie Symmetry method\cite{Sahoo2017}, Riccati equation method\cite{Malik2023}, modified extended direct algebraic method\cite{Seadawy2017}, new extended generalized method\cite{Tripathy2024,Zayed2020},Unified generalized  kudryashov methhod\cite{Darwish2007}, ($\frac{G^{\prime}}{G^{\prime}+G+A}$) method\cite{Tripathy2020}, Extended hyperbolic function method\cite{REHMAN2022105802}, 
	 He's homotopy perturbation method\cite{Saberi-Nadjafi2009} etc.\\
	
There is another interest growing among many researchers in describing the qualitative behavior of the system using dynamical analysis, because it reveals the stability, evolution, and structural behavior of the obtained solution. In this analysis, the phase portrait of the system illustrates the trajectories of the achieved results, which demonstrate the behavior of the system, the sensitive analysis that displays the sensitivity of a system by varying the initial condition, and the bifurcation of a system that describes how the stability changes of a system vary by varying the parameter \cite{Wazwaz2007,Tripathy2024,Zayed2020}.
	
	In general, most nonlinear equations begin with the KdV equation\cite{AMW,Lakshmanan2003} which delineates waves in shallow water.\\
	\begin{equation}
		u_t+6u u_x+u_{xxx}=0
	\end{equation}\\
	But this equation is limited to one-dimensional and low-order dispersion; as a result it fails to describe the complexity of the system.\\
	
	In 1970, Kadomtsev and Petviashvili added a term to the KdV equation \cite{Lakshmanan2003} that helps to move the wave primarily in $x$ but evolve slightly in $y$,  which is given below:
	\begin{equation}
		(u_t+6u u_x+u_{xxx})_x+ \alpha u_{yy}=0
	\end{equation}\\
	Later, in 1976, Caudrey, Dodd, and Gibbon identified a specific hierarchy of the KDV equation; they established a fifth-order and seventh-order equation named the Caudrey-Dodd-Gibbon equation. In \cite{Abdollahzadeh2010} M.Abdollahzadeh, M. Hosseini, M.Ghanbarpour and H.Shirvani investigate the fifth-order  Caudrey-Dodd-Gibbon equation and in \cite{Sharma2020} Ankita Sharma and Rajan Arora investigate the seventh-order  Caudrey-Dodd-Gibbon equation, which are respectively given as:
		\begin{equation}
		u_t + 30\,u_x u_{xx} + 30\,u\,u_{xxx} + 180\,u^2 u_x + u_{xxxxx} = 0
	\end{equation}
	and
	\begin{equation}
		u_t 
		+ 420 u^{3} u_x 
		+ 210 u^{2} u_{3x} 
		+ 420 u u_x u_{2x} 
		+ 28 u u_{5x} 
		+ 28 u_x u_{4x} 
		+ 70 u_{2x} u_{3x} 
		+ u_{7x}
	\end{equation}\\
	In this paper, the core focus of this work is to investigate the (2 + 1)-dimensional seventh-order Caudrey-Dodd-Gibbon-KP (sCDG-KP) equation, which originated by merging the seventh-order CGP equation and the KP equation.\\
	\begin{equation}\label{1.5}
		\begin{aligned}
		\left(
		u_t 
		+ 420 u^{3} u_x 
		+ 210 u^{2} u_{3x} 
		+ 420 u u_x u_{2x} 
		+ 28 u u_{5x} 
		+ 28 u_x u_{4x} 
		+ 70 u_{2x} u_{3x} 
		+ u_{7x}
		\right)_x \\
		+ \alpha u_{2y} = 0.
		\end{aligned}
	\end{equation}
	where, $\alpha = \pm 1$.
	
	In\cite{Qin2024}, Mengyao Qin, Yunhu Wang, and Manwai Yuen establish its Lie symmetry and find its exact solution using a unified algebric method.\\
	
	The whole work is organized into four distinct sections, such as
	Section 2 presents the description of ($\frac{G^{\prime}}{G^{\prime}+G+A}$) method and utilization of the given method to obtain the exact solution of Eq.(\ref{1.5}). Section 3 demonstrates the dynamical analysis of Eq.(\ref{1.5}) by converting Eq.(\ref{1.5}) into a 2D system of equations and portrays the sensitivity of the governing system to the variation of the parameter and initial condition. Lastly, section 4 describes the conclusion of the whole work.

	\section{Description of ($\frac{G^{\prime}}{G^{\prime}+G+A}$) method:}
	Let us consider the following nonlinear PDE:
	\begin{equation}\label{2.1}
		F(u,u_{x},u_{y},u_{z},u_{t},u_{xx},u_{yy},u_{zz},u_{xy}u_{tt},.......)=0,
	\end{equation}
	where $F$ is a polynomial of the unknown function $u=u(x,y,z,t)$ and $u_{x},u_{y},u_{z},u_{t}$..... represent the partial derivatives of the function $u(x,y,z,t)$. Here, $u$ represent as the dependent variable and $x,y,z,$ and $t$ are the independent variables.\\
	\textbf{Step:1:-}
	Let us, define the traveling wave transformation as:
	\begin{equation}\label{2.2}
		u(x,y,t)=u(\xi),  
	\end{equation}
	
	where $\xi=x+\omega y-\sigma t$ and $ \omega, \sigma $, are non-zero constants.
	
	The Eq.(\ref{2.1}) is now converted into the following ODE via the Eq.(\ref{2.2}) and the wave transform  $\xi=x+\omega y-\sigma t$.
	
	\begin{equation}\label{2.3}
		Q(u^{\prime},u^{\prime \prime},u^{\prime\prime\prime},......)=0.
	\end{equation}\\
	\textbf{Step:2:-}
	Let us assume that Eq.(\ref{2.3}) has an exact traveling wave solution is in the following form:
	\begin{equation}\label{2.4}
		u({\xi})=\sum_{j=0}^{N} a_j \left( \frac{G_0}{G_0 + G + A} \right)^j
	\end{equation}
	Where,
	$a_{0},a_{1},a_{2},...a_{N}$ are arbitrary constants, and the value of $N$ can be determined via the homogeneous balancing method.\\Here, $G(\xi)$ is the solution of the following second-order linear Ordinary differential equation:
	\begin{equation}\label{2.5}
		G'' + B G' + C G + A C = 0.
	\end{equation}\\
	Where, $B, C$, and $A$ are the real constants.

	\textbf{Step:3:-} Now, by substituting Eq. (\ref{2.4}) along with Eq. (\ref{2.5}) into Eq. (\ref{2.3}), a polynomial having the power of ($\frac{G^{\prime}}{G^{\prime}+G+A}$) is formulated. By collecting each power of ($\frac{G^{\prime}}{G^{\prime}+G+A}$) and setting its coefficients equal to zero, a system of equations is generated. After solving this system of equations, we obtain the value of  $a_{0}, a_{1}, a_{2},...a_{N}$,  $B, C$, and $A$. We can obtain the exact solution of NLPDE by using the value of $a_{0}, a_{1}, a_{2},...a_{N}$,  $B, C$, and $A$\cite{Tripathy2020}.
	\subsection{Traveling Wave Solution:}
	In this section, we examine the solution of Eq.(\ref{1.5}) using ($\frac{G^{\prime}}{G^{\prime}+G+A}$) method. let us consider the following (2 + 1)-dimensional seventh-order Caudrey-Dodd-Gibbon-KP (sCDG-KP) equation as:\\
	$
	\left(
	u_t 
	+ 420 u^{3} u_x 
	+ 210 u^{2} u_{3x} 
	+ 420 u u_x u_{2x} 
	+ 28 u u_{5x} 
	+ 28 u_x u_{4x} 
	+ 70 u_{2x} u_{3x} 
	+ u_{7x}
	\right)_x
	+ \alpha u_{2y} = 0,
	\quad \alpha = \pm 1.
	$\\
	Substituting the wave transform  $\xi=x+\omega y+qz-\sigma t$ and Eq.(\ref{2.2}) into Eq.(\ref{1.5}), then we obtain the following ODE as:
	\begin{equation}\label{2.6}
		\begin{aligned}
			&1260 u^2 (u')^2
			+ 420 u^3 u''
			+ 420 (u')^2 u''
			+ 420 u (u'')^2
			+ 840 u u' u^{(3)} \\
			&\quad
			+ 70 (u^{(3)})^2
			+ 210 u^2 u^{(4)}
			+ 98 u'' u^{(4)}
			+ 56 u' u^{(5)} \\
			&\quad
			+ 28 u u^{(6)}
			+ u^{(8)}
			+ (\alpha \omega^{2} - \sigma) u'' = 0.
		\end{aligned}
	\end{equation}
	Integrating twice of the above equation w.r.t $\xi$ and simplifying, gives:
	\begin{equation}\label{2.7}
		u^{(6)}+105 u^{(4)}+210 u^2 u''+28 u u^{(4)}+35 u''^{2}+(\alpha \omega^{2} - \sigma) u'=0
	\end{equation}
	Using homogeneous balancing principle on the Eq.(\ref{2.7}) we find $N=2$, then according to ($\frac{G^{\prime}}{G^{\prime}+G+A}$) method we get:
	\begin{equation}\label{2.8}
		u({\xi})= a_0+a_1 \left( \frac{G_0}{G_0 + G + A} \right)+a_2  \left( \frac{G_0}{G_0 + G + A} \right)
	\end{equation}
	By swapping Eq.(\ref{2.8}) and Eq.(\ref{2.5}) into Eq.(\ref{2.7}) a system of equation is formulated by collecting each power of ($\frac{G^{\prime}}{G^{\prime}+G+A}$) and setting it's coefficients is equal to zero and after solving this, we derive the following result as:
	\begin{equation}\label{2.9}
		\begin{aligned}
			a_{0} &= -2 \left( C - BC + C^{2} \right), \\
			a_{1} &= 2 \left( -B + B^{2} + 2C - 3BC + 2C^{2} \right), \\
			a_{2} &= -2 \left( -1 + B - C \right)^{2}, \\
			\sigma &= B^{6} - 12 B^{4} C + 48 B^{2} C^{2} - 64 C^{3} + \alpha \omega^{2}.
		\end{aligned}
	\end{equation}\\
	Upon putting these value into Eq.(\ref{2.8}) and Eq.(\ref{2.2}), we obtain the exact solution of (2 + 1)-dimensional seventh-order Caudrey-Dodd-Gibbon-KP (sCDG-KP) equation, which is represented in Eq.(\ref{1.5}):\\
	\textbf{\underline{Case:1}}
	When, $\Gamma=B^2-4C>0$
	\begin{equation}\label{2.10}
		\begin{aligned}
			u(\xi) =\;& -2 \left(C - BC + C^{2}\right) \\
			&+ \frac{
				2 \left(-B + B^{2} + 2C - 3BC + 2C^{2}\right)
				\left[(B + \sqrt{\Gamma}) c_{1}
				+ e^{\sqrt{\Gamma}\xi} (B - \sqrt{\Gamma}) c_{2}\right]
			}{
				(-2 + B + \sqrt{\Gamma}) c_{1}
				+ e^{\sqrt{\Gamma}\xi} (-2 + B - \sqrt{\Gamma}) c_{2}
			} \\
			&- \frac{
				2 (-1 + B - C)^{2}
				\left[(B + \sqrt{\Gamma}) c_{1}
				+ e^{\sqrt{\Gamma}\xi} (B - \sqrt{\Gamma}) c_{2}\right]^{2}
			}{
				\left[(-2 + B + \sqrt{\Gamma}) c_{1}
				+ e^{\sqrt{\Gamma}\xi} (-2 + B - \sqrt{\Gamma}) c_{2}\right]^{2}
			}.
		\end{aligned}
	\end{equation}
	\textbf{\underline{Case:2}}
	When, $\Gamma=B^2-4C<0$
	
	\begin{equation}\label{2.11}
		\begin{aligned}
			u(\xi) =\;& -2 \left( C - BC + C^{2} \right) + 
			2 \left( -B + B^{2} + 2C - 3BC + 2C^{2} \right) \\
			& \frac{
				\left[
				\sin\!\left( \frac{\sqrt{-\Gamma}\,\xi}{2} \right)
				\left( \sqrt{-\Gamma}\, C_{1} + B C_{2} \right)
				+
				\cos\!\left( \frac{\sqrt{-\Gamma}\,\xi}{2} \right)
				\left( B C_{1} - \sqrt{-\Gamma}\, C_{2} \right)
				\right]
			}{
				\sin\!\left( \frac{\sqrt{-\Gamma}\,\xi}{2} \right)
				\left( \sqrt{-\Gamma}\, C_{1} + (B-2) C_{2} \right)
				+
				\cos\!\left( \frac{\sqrt{-\Gamma}\,\xi}{2} \right)
				\left( (B-2) C_{1} - \sqrt{-\Gamma}\, C_{2} \right)
			} \\
			&- \frac{
				2 (-1 + B - C)^{2}
				\left[
				\sin\!\left( \frac{\sqrt{-\Gamma}\,\xi}{2} \right)
				\left( \sqrt{-\Gamma}\, C_{1} + B C_{2} \right)
				+
				\cos\!\left( \frac{\sqrt{-\Gamma}\,\xi}{2} \right)
				\left( B C_{1} - \sqrt{-\Gamma}\, C_{2} \right)
				\right]^{2}
			}{
				\left[
				\sin\!\left( \frac{\sqrt{-\Gamma}\,\xi}{2} \right)
				\left( \sqrt{-\Gamma}\, C_{1} + (B-2) C_{2} \right)
				+
				\cos\!\left( \frac{\sqrt{-\Gamma}\,\xi}{2} \right)
				\left( (B-2) C_{1} - \sqrt{-\Gamma}\, C_{2} \right)
				\right]^{2}
			}.
		\end{aligned}
	\end{equation}
	
\begin{figure}[htbp]
		\centering
		
		\begin{subfigure}{0.55\textwidth}
			\centering
			\includegraphics[width=.75\linewidth]{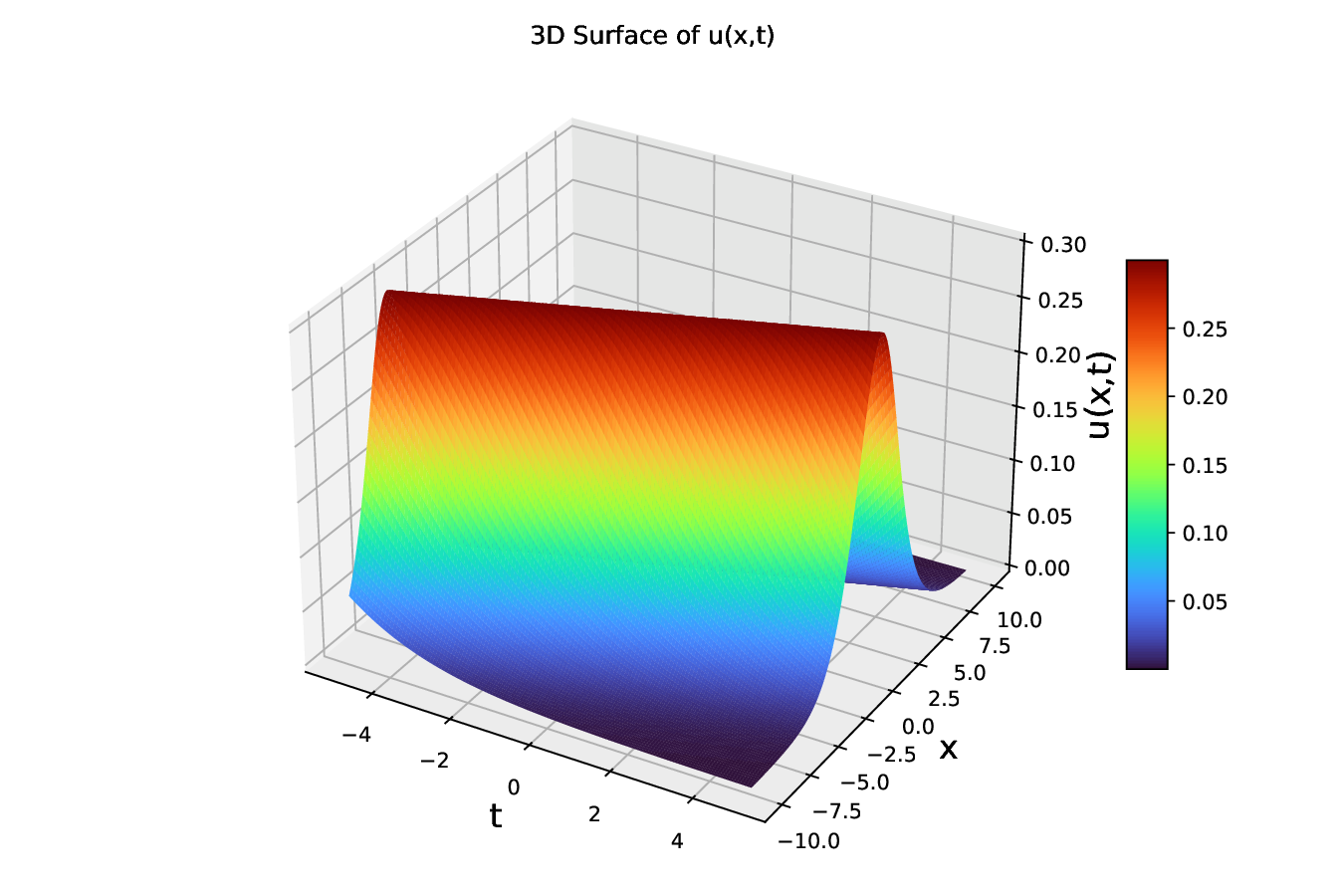}
			\caption{}
			\label{fig1a}
		\end{subfigure}
		\hfill
		\begin{subfigure}{0.45\textwidth}
			\centering
			\includegraphics[width=.75\linewidth]{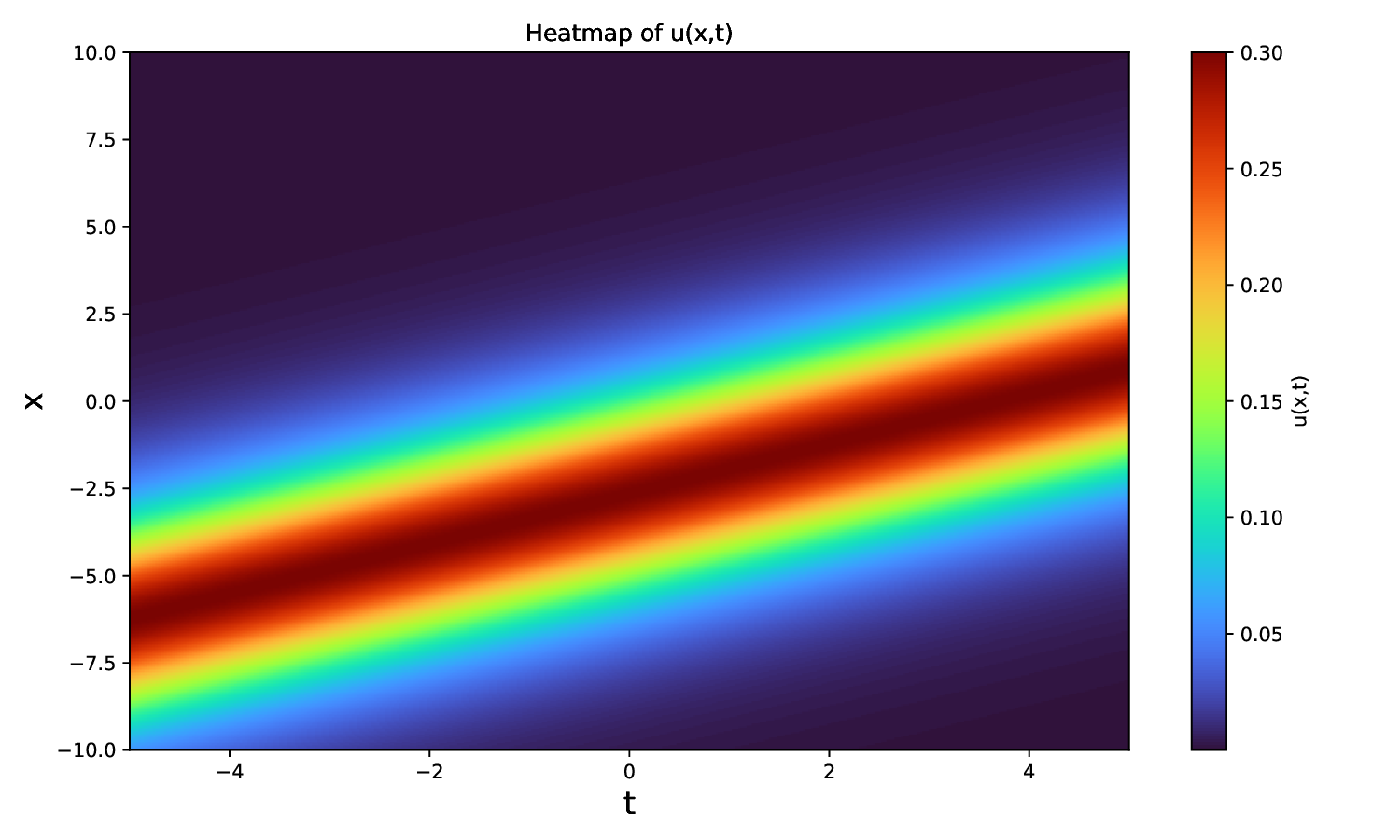}
			\caption{}
			\label{fig1b}
		\end{subfigure}
		\begin{subfigure}{0.45\textwidth}
			\centering
			\includegraphics[width=.75\linewidth]{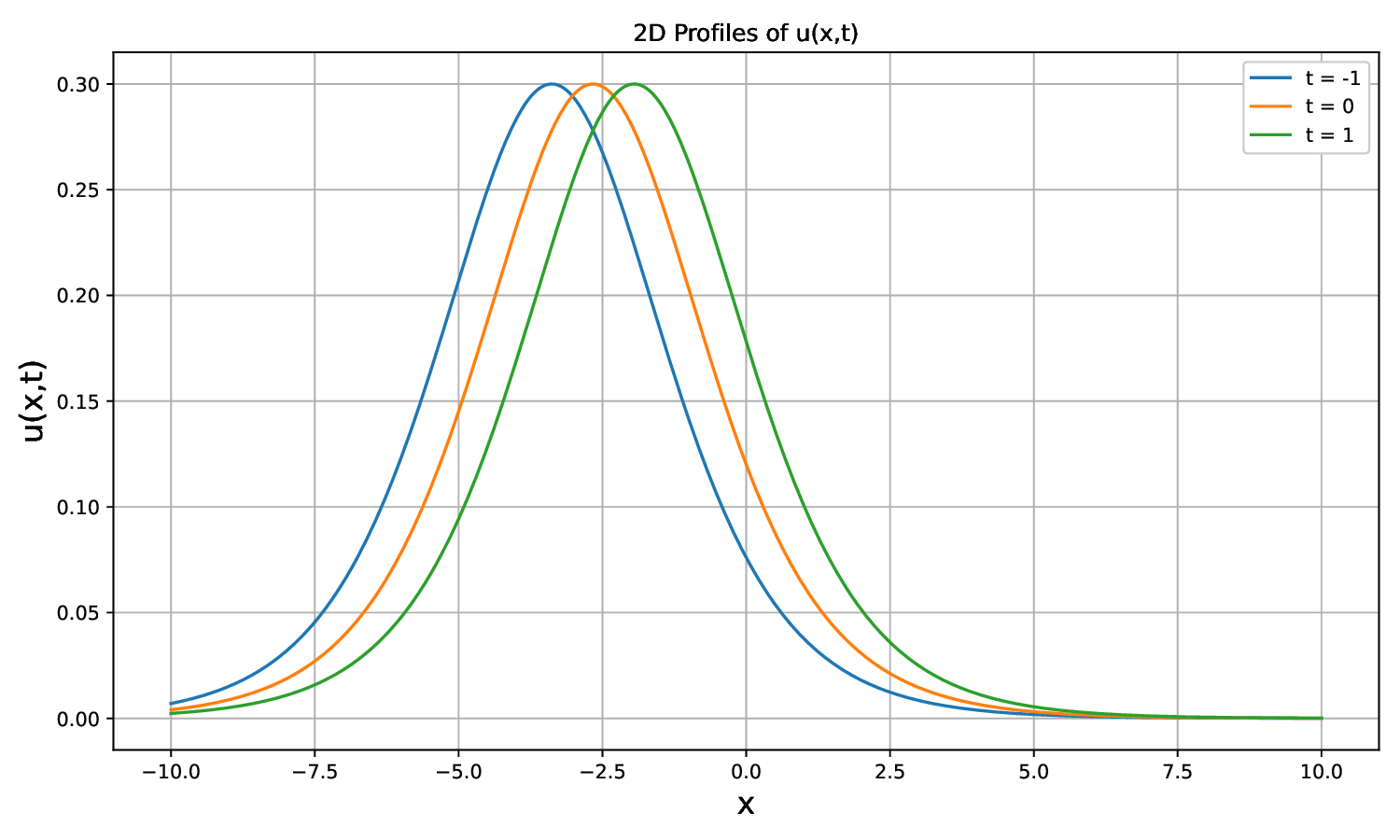}
			\caption{}
			\label{fig1c}
		\end{subfigure}
		
		\caption{Graphical representation of Eq,\ref{2.10} with $B=1$, $C=1$, $\alpha=0.5$, $\omega=1$, $C_1=1$, $C_2=1$. }
		\label{fig1:main}
	\end{figure}

\begin{figure}[htbp]
		\centering
		
		\begin{subfigure}{0.55\textwidth}
			\centering
			\includegraphics[width=.75\linewidth]{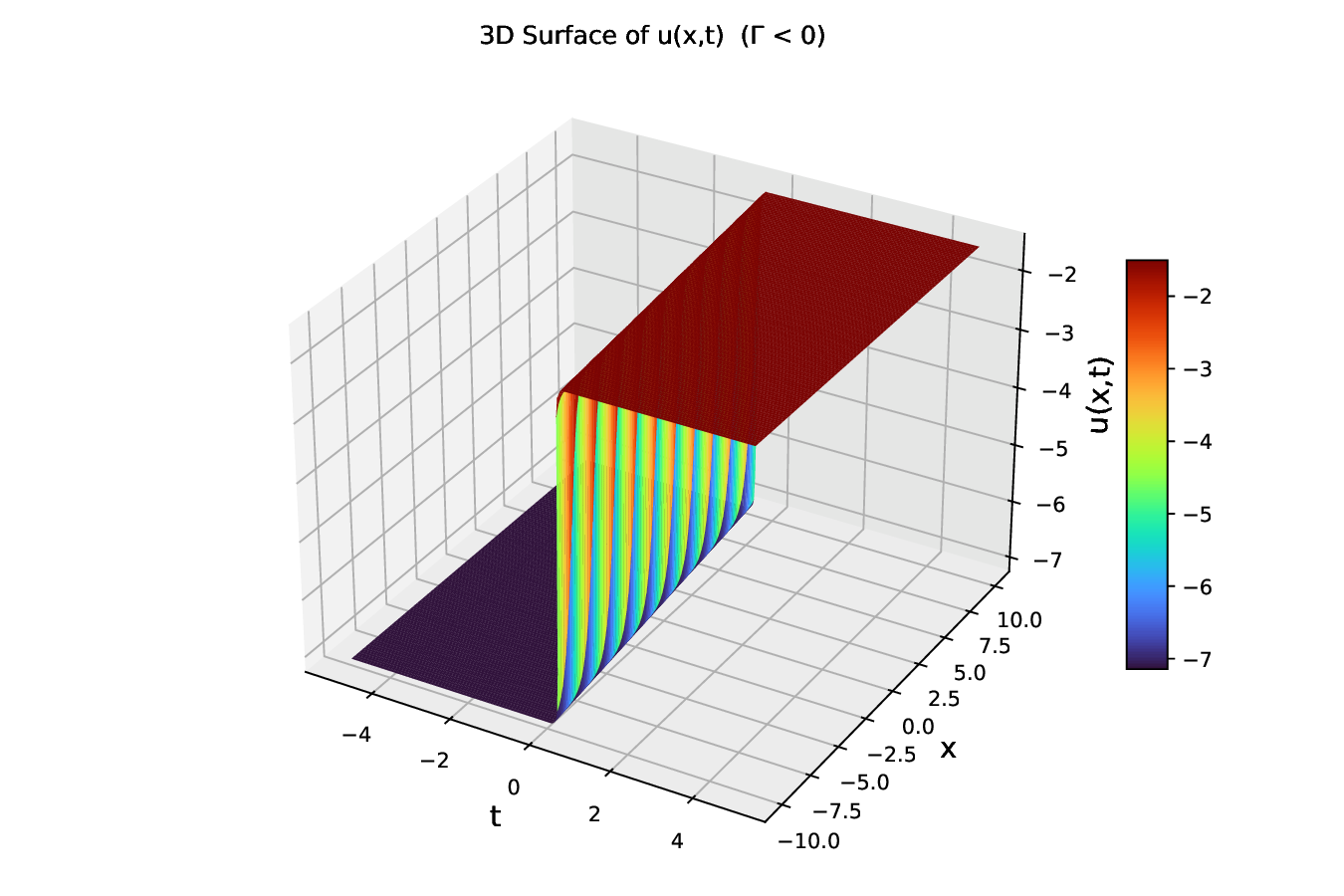}
			\caption{}
			\label{fig2a}
		\end{subfigure}
		\hfill
		\begin{subfigure}{0.45\textwidth}
			\centering
			\includegraphics[width=.75\linewidth]{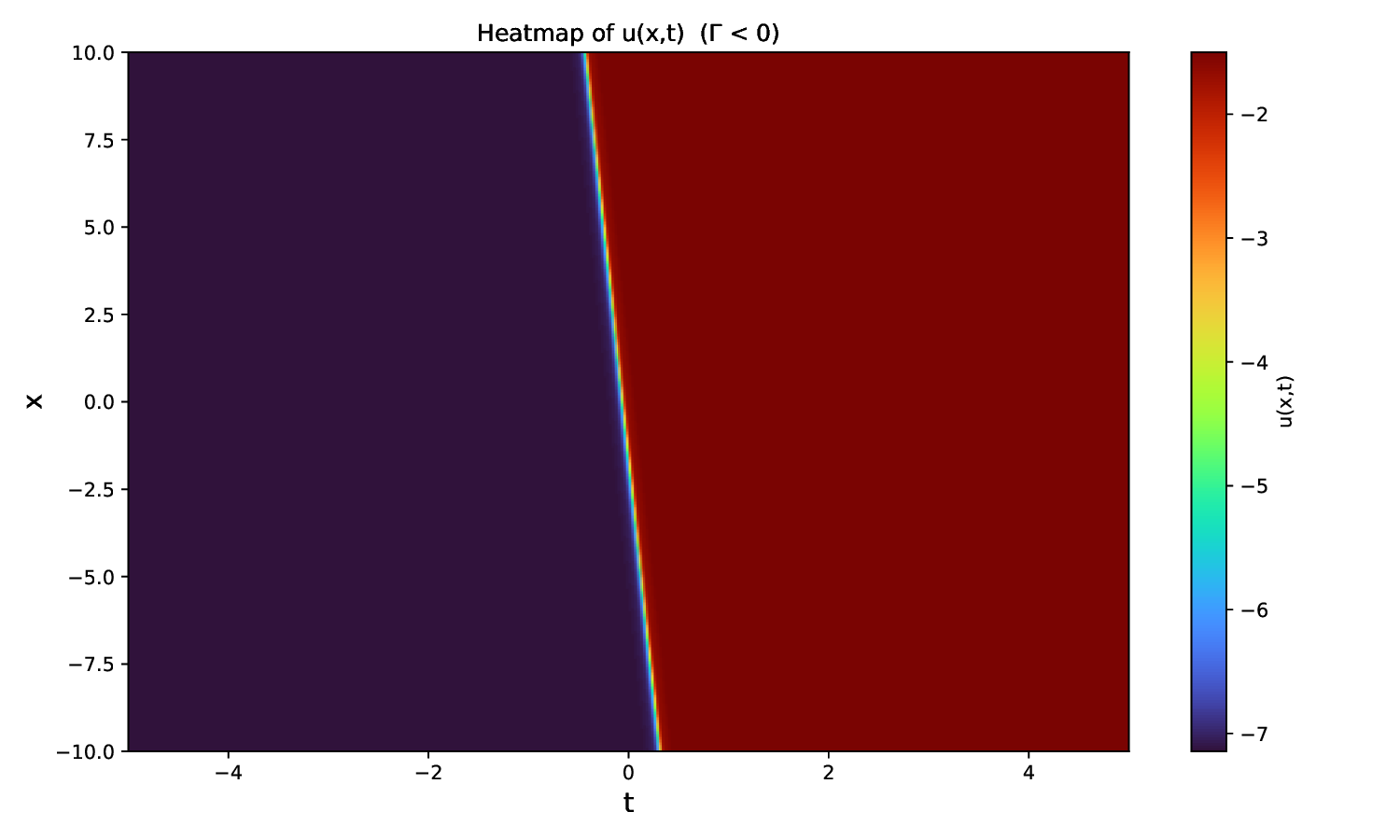}
			\caption{}
			\label{fig2b}
		\end{subfigure}
		\begin{subfigure}{0.45\textwidth}
			\centering
			\includegraphics[width=.75\linewidth]{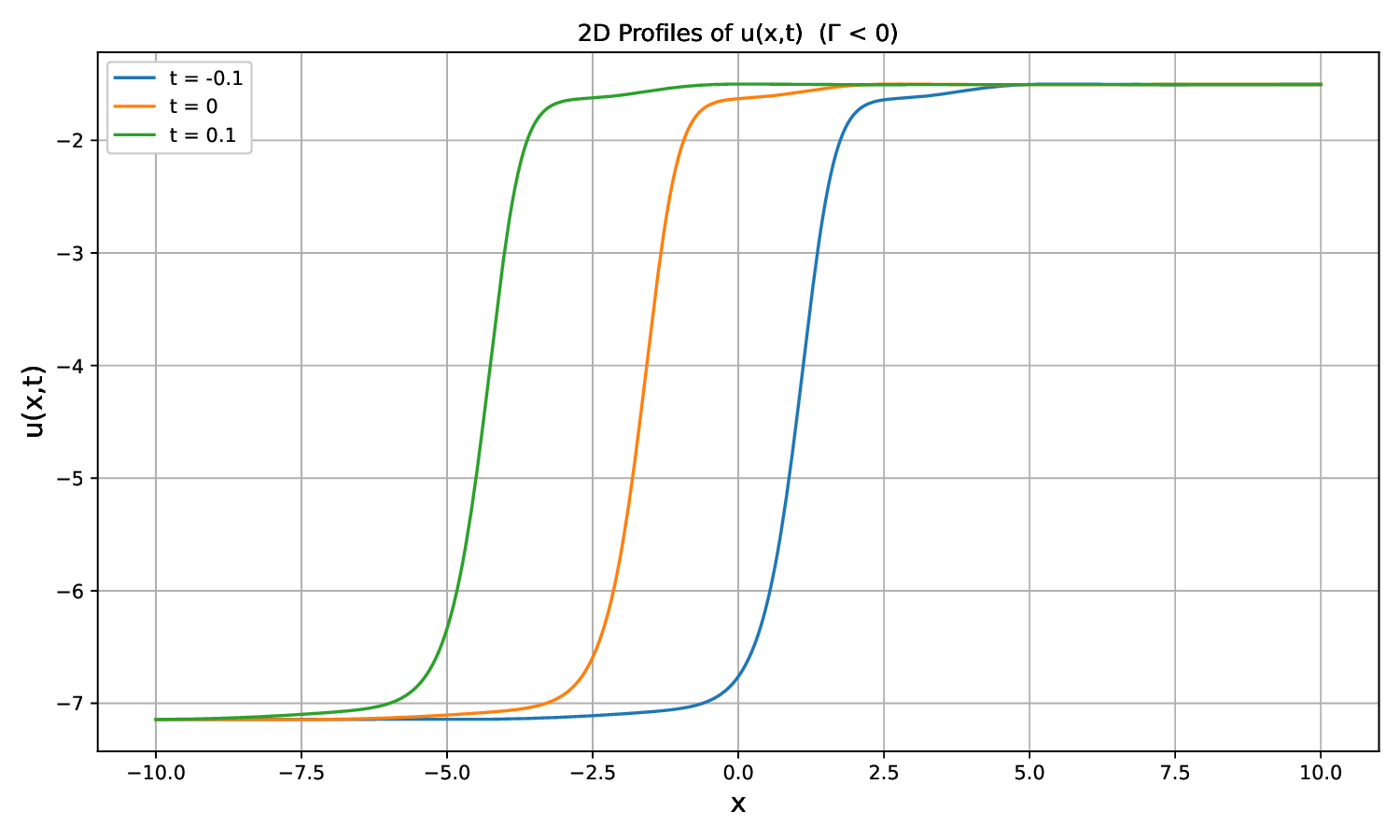}
			\caption{}
			\label{fig2c}
		\end{subfigure}
		
		\caption{Graphical representation of Eq,\ref{2.11} with $B=1$, $C=1$, $\alpha=0.5$, $\omega=\sqrt{-1}$, $C_1=\sqrt{-2.1}$, $C_2=-1$.}
		\label{fig2:main}
	\end{figure}

	Fig.(\ref{fig1:main}) and fig.(\ref{fig2:main}) illustrate the graphical visualization of Eq.(\ref{2.10}) and (\ref{2.11}) respectively.
	
	\section{\textbf{Dynamical Analysis:-}}
	This section demonstrates the qualitative behavior of Eq.(\ref{1.5}) with the help of different tools such as the phase space analysis, bifurcation analysis, sensitive analysis and indicators of chaos.
	Now, we have to convert the Eq.(\ref{1.5}) into two systems of equations by the using following procedure.\\
	Let us consider Eq.(\ref{2.7})
	
	$u^{(vi)} + 105u^{(iv)} + 210u^{2}u^{\prime\prime} + 28u u^{(iv)}+ \left(u''\right)^2
	+ ( \alpha \omega^2  - \sigma) u=0$.
	
	With the help of the polynomial trail method, we can take a trail equation as: 
	
	\begin{equation}\label{3.1}
		u^{\prime\prime}=\sum_{k=0}^{N}p_{k}u^{k}.
	\end{equation}
	
	By homogeneous balancing method we know that $N$=2, so Eq.(\ref{3.1}) becomes
	
	\begin{equation}\label{3.2}
		u^{\prime\prime}= p_{0}+ p_{1}u+p_{2}u^{2}.
	\end{equation} 
	
	Integrating both side With respect to $\xi$ of the above  Eq.(\ref{3.2}) we get,
	
	\begin{equation}
		(u^{\prime})^{2}=2p_{0}u+p_{1}u^{2}+\frac{2p_{2}}{3}u^{3}+d_{0}.
	\end{equation}
	
	\begin{equation}\label{3.4}
		\Rightarrow	(u^{\prime})=\sqrt{2p_{0}u+p_{1}u^{2}+\frac{2p_{2}}{3}u^{3}}+d_{0}.
	\end{equation}

	Where $ d_{0}$ is an integration constant.
	Now, by putting  Eq.(\ref{3.2}) along with  Eq.(\ref{3.4}) into Eq.(\ref{2.7}), a system of equations is generated by collecting each power of $u$ such as:
	
	\begin{equation}\label{3.5}
		\begin{aligned}
			&p_{0} p_{1}^2 + 10 d_0 p_{1} p_{2} + p_{0}^2 (35 + 6 p_{2}) = 0, \\
			&-\sigma + \alpha \omega^2 + p_{1}^3 
			+ 4 d_0 p_{2} (14 + 5 p_{2}) 
			+ 2 p_{0} p_{1} (49 + 18 p_{2}) = 0, \\
			&21 p_{1}^2 (3 + p_{2}) 
			+ 14 p_{0} (15 + 17 p_{2} + 4 p_{2}^2) = 0, \\
			&\frac{70}{3} p_{1} (9 + 9 p_{2} + 2 p_{2}^2) = 0, \\
			&\frac{35}{3} (9 + 18 p_{2} + 11 p_{2}^2 + 2 p_{2}^3) = 0.
		\end{aligned}
	\end{equation}
	After solving the above system of Eq.(\ref{3.5}) we get:
	\begin{equation}\label{3.6}
		\begin{aligned}
			p_{0} &= \frac{3(-\sigma + \alpha \omega^{2})^{2/3}}{2 \times 89^{2/3}}, \\
			p_{1} &= \frac{(-\sigma + \alpha \omega^{2})^{1/3}}{89^{1/3}}, \\
			p_{2} &= -\frac{3}{2}, \\
			d_{0} &= -\frac{4}{89}(\sigma - \alpha \omega^{2}).
		\end{aligned}
	\end{equation}
	Here we can clearly see that $p_2$ is always negative.
	\subsection{Bifurcation analysis:}
	Phase portrait analysis presents the portrait of the behavior of a dynamical system by analyzing the trajectories of its solution. In this work, we apply bifurcation analysis, which provides insight into the qualitative behavior of the Eq.(\ref{1.5}) that shows how a small variation of the parameter causes a significant change in the dynamic system\cite{Almheidat2025,Nasreen2024}. \\
	Let us take a new notation such as $u'=A$. Then Eq(\ref{3.2}) can be transformed in 2D dynamical system as:
	\begin{equation}\label{3.7}
		\begin{cases}
			\dfrac{du}{d\xi} = A, \\[6pt]
			\dfrac{dA}{d\xi} = p_{2}u^{2} + p_{1}u + p_{0}.  \\[6pt]
		\end{cases}
	\end{equation}
	The above system satisfies the Hamiltonian characteristics and possesses the following:
	\begin{equation}\label{3.8}
		\frac{A^2}{2}-(p_2 \frac{u^3}{3}+p_1 \frac{u^2}{2}+p_0 u)
	\end{equation}
	The equilibrium points of Eq.(\ref{3.7}) are:\\
	
\begin{figure}[htbp]
		\centering
		
		\begin{subfigure}{0.45\textwidth}
			\centering
			\includegraphics[width=.5\linewidth]{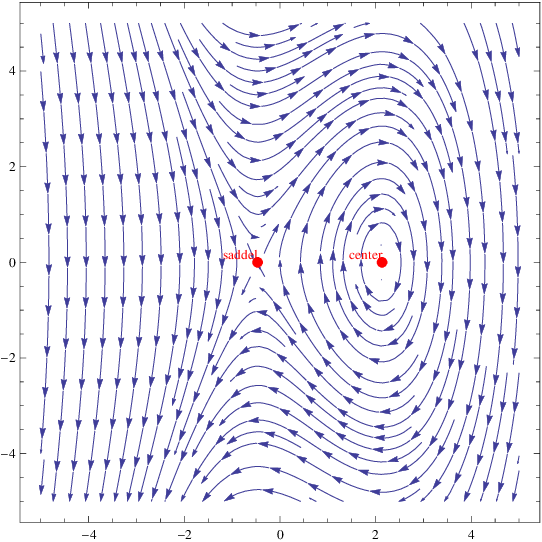}
			\caption{ $p_{0}, p_{1} > 0$}
			\label{fig3a}
		\end{subfigure}
		\hfill
		\begin{subfigure}{0.45\textwidth}
			\centering
			\includegraphics[width=.5\linewidth]{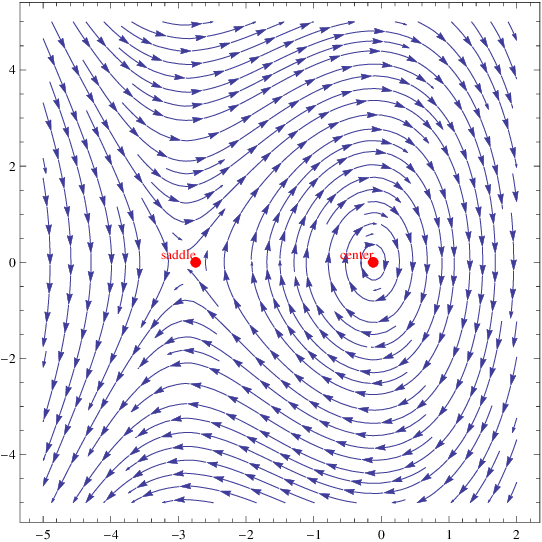}
			\caption{$p_{0}, p_{1} < 0$}
			\label{fig3b}
		\end{subfigure}
		\hfill
		\begin{subfigure}{0.45\textwidth}
			\centering
			\includegraphics[width=.5\linewidth]{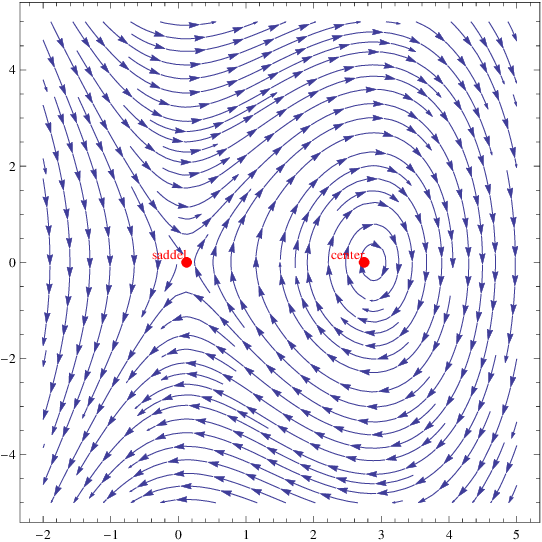}
			\caption{$p_{0}>, p_{1} < 0$}
			\label{fig3c}
		\end{subfigure}
		\hfill
		\begin{subfigure}{0.45\textwidth}
			\centering
			\includegraphics[width=.5\linewidth]{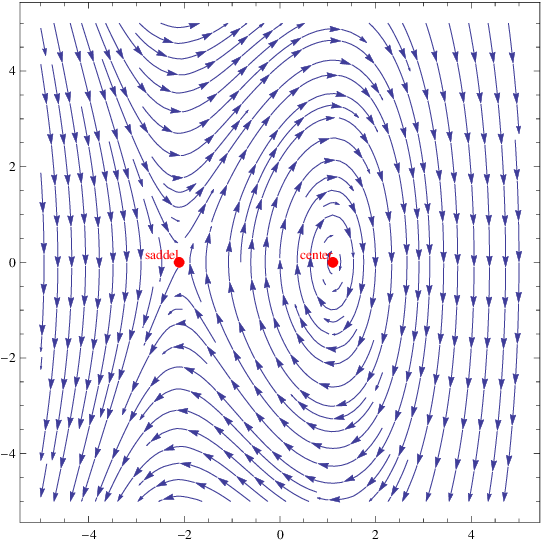}
			\caption{$p_{0}<, p_{1} > 0$}
			\label{fig3d}
		\end{subfigure}

		\caption{Bifurcation portrait of system(\ref{3.7}) with different condition. }
		\label{fig:main}
	\end{figure}

	$
	E_{1} = (\frac{-p_{1} + \sqrt{p_{1}^{2} - 4 p_{2} p_{0}}}{2 p_{2}},0)$ and
	$
	E_{2} = (\frac{-p_{1} - \sqrt{p_{1}^{2} - 4 p_{2} p_{0}}}{2 p_{2}},0).$\\
	And the determinant of the Jacobian matrix of Eq.(\ref{3.7}) is given as:
	\begin{equation}
		J(u,A)=-s_{1}-2s_{2}u.
	\end{equation}
	Now we examine all possible cases of equilibrium points of Eq.(\ref{3.7}), since $s_2$ is always negative, so four different cases arise as:\\
	(1) $p_{0},p_1<0$.\\
	(2)$p_{0},p_1>0$.\\
	(3)$p_{0}>0,p_1<0$.\\
	(4))$p_{0}<0,p_1>0.$\\
	\textbf{Case:1}\\
	when $p_{0},p_1>0$, the two equilibrium points becomes $E_{1}=(-0.468, 0)$ and $ E_{2}=(2.135, 0)$ by taking $p_{0}=1.5$ and $p_1=2.5$, in fig.(\ref{fig3a}) we can clearly see that $E_{1}$ possess saddle point and $ E_{2}$ posses center point.\\
	\textbf{Case:2:-}\\
	when and $p_0,p_1<0$ the two equilibrium points becomes $E_{1}=(-0.121, 0)$ and $ E_{2}=(2.745, 0)$ by taking  $p_{0}=-0.5$ and $p_1=-4.5$ , in fig.(\ref{fig3b}) we can see that  $M_{1}$ possess center point and $ M_{2}$ posses saddle point.\\
	\textbf{Case:-3}\\
	when $p_{0}>0$ and,$p_1<0$ the two equilibrium points becomes $E_{1}=(-2.107, 0)$ and $ E_{2}=(1.107, 0)$ by taking $p_{0}=3.5$ and $p_1=-1.5$,in fig.(\ref{fig3c}) we can see that $M_{1}$ possess center point and $ M_{2}$ posses saddle point.\\
	\textbf{Case:-4}\\
	when $p_{0}<0$ and,$p_1>0$ the two equilibrium points becomes $E_{1}=(0.121, 0)$ and $ E_{2}=(2.745, 0)$ by taking $p_{0}=-0.5$ and $p_1=4.5$, in fig.(\ref{fig3d}) we can see that  $M_{1}$ possess center point and $ M_{2}$ posses saddle point.\\

	\subsection{Chaotic analysis:}
	In this section, we will explore the dynamical system \ref{3.7} by adding a perturbation term to it. The following perturbed term is given as\cite{Nasreen2024}:
	\begin{equation}\label{3.10}
		\begin{cases}
			\dfrac{du}{d\xi} = A, \\[6pt]
			\dfrac{dA}{d\xi} = p_{2}u^{2} + p_{1}u + p_{0}+Z_0 \cos(\nu \xi).  \\[6pt]
		\end{cases}
	\end{equation}
	Here, the perturbed term is denoted as $Z_0$ and $\nu$ represent the frequency.\\
	
\begin{figure}[htbp]
		\centering
		
		\begin{subfigure}{0.45\textwidth}
			\centering
			\includegraphics[width=.5\linewidth]{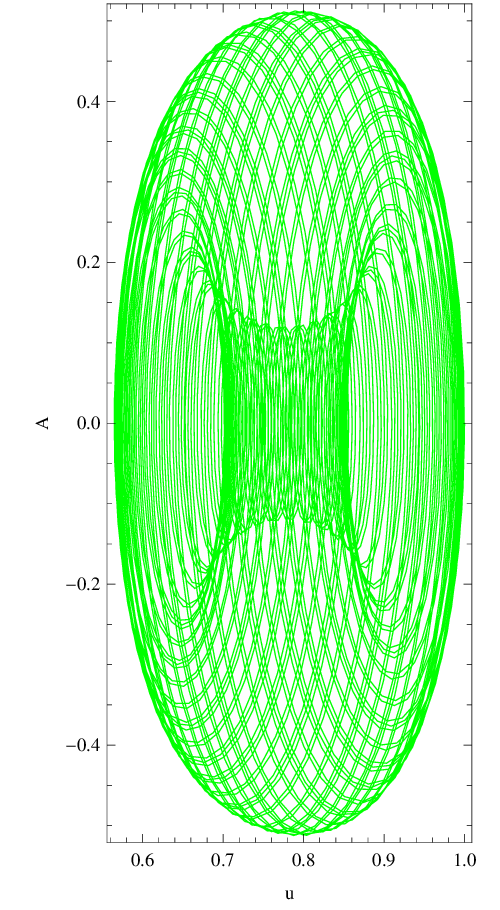}
			\caption{2D phase portrait}
			\label{fig4a}
		\end{subfigure}
		\hfill
		\begin{subfigure}{0.45\textwidth}
			\centering
			\includegraphics[width=.75\linewidth]{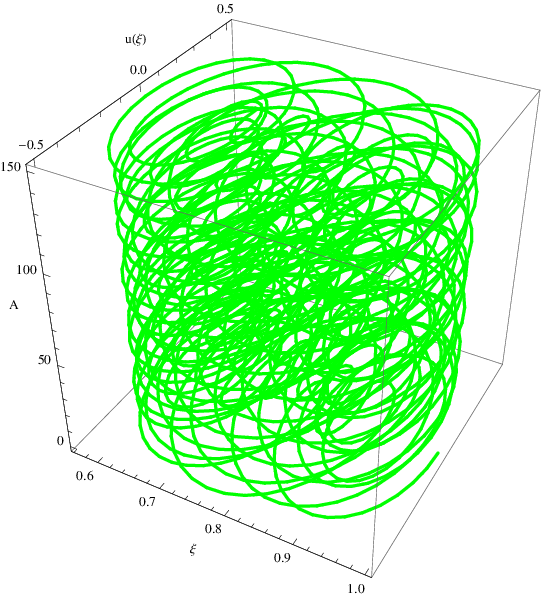}
			\caption{3D phase portrait}
			\label{fig4b}
		\end{subfigure}
		\begin{subfigure}{0.45\textwidth}
			\centering
			\includegraphics[width=.75\linewidth]{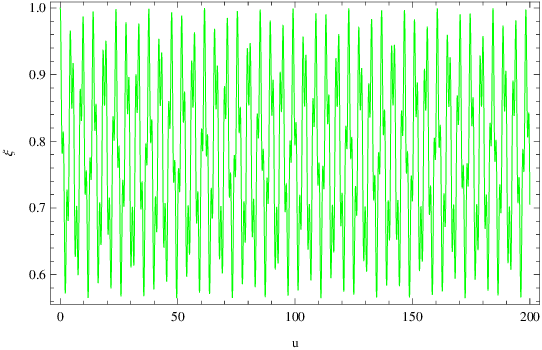}
			\caption{Time series }
			\label{fig4c}
		\end{subfigure}
		\caption{Phase portrait of system(\ref{3.10}) with $\sigma = -1.1$, $\alpha = 1$, $\omega = 3.9$, $z_0 =
			-1.3$, $\nu = -4.5$. }
		\label{fig4:main}
	\end{figure}

\begin{figure}[htbp]
		\centering
		
		\begin{subfigure}{0.35\textwidth}
			\centering
			\includegraphics[width=.55\linewidth]{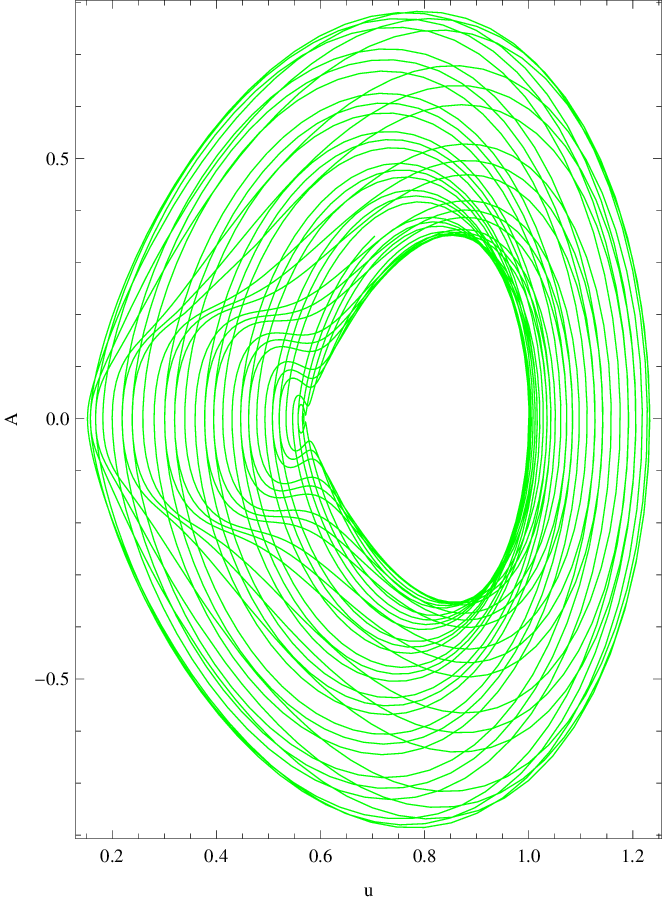}
			\caption{2D phase portrait}
			\label{fig5a}
		\end{subfigure}
		\hfill
		\begin{subfigure}{0.45\textwidth}
			\centering
			\includegraphics[width=.55\linewidth]{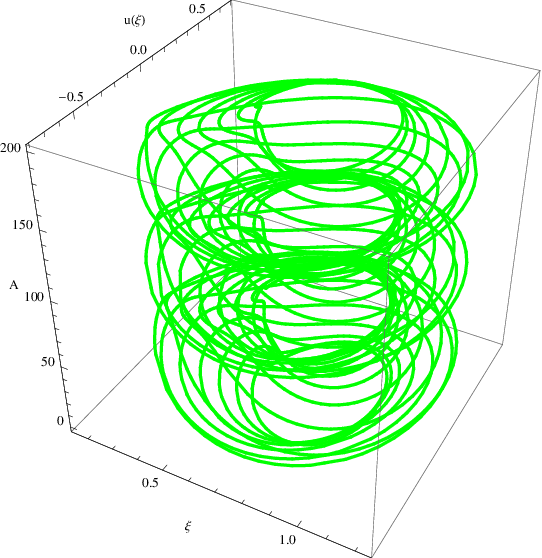}
			\caption{3D phase portrait}
			\label{fig5b}
		\end{subfigure}
		\begin{subfigure}{0.45\textwidth}
			\centering
			\includegraphics[width=.85\linewidth]{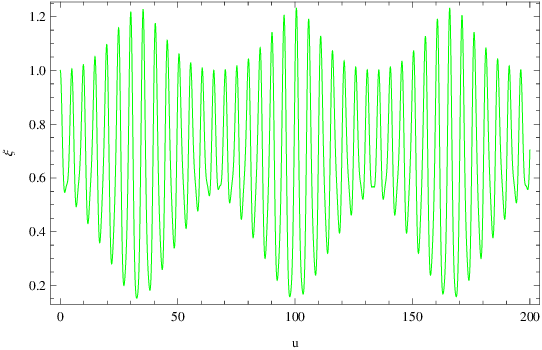}
			\caption{Time series }
			\label{fig5c}
		\end{subfigure}
		\caption{Phase portrait of system(\ref{3.10}) wit $\sigma = -1.05$, $\alpha = 1$, $\omega = -3.5$, $z_0 =
			-0.3$, $\nu = -2.5$. }
		\label{fig5:main}
	\end{figure}

\begin{figure}[htbp]
		\centering
		
		\begin{subfigure}{0.35\textwidth}
			\centering
			\includegraphics[width=.55\linewidth]{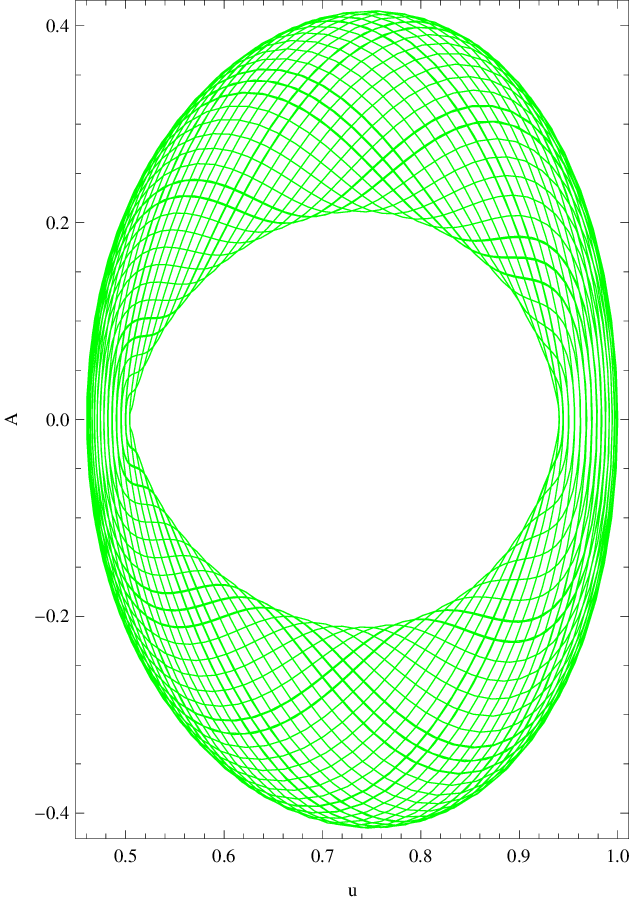}
			\caption{2D phase portrait}
			\label{fig6a}
		\end{subfigure}
		\hfill
		\begin{subfigure}{0.45\textwidth}
			\centering
			\includegraphics[width=.55\linewidth]{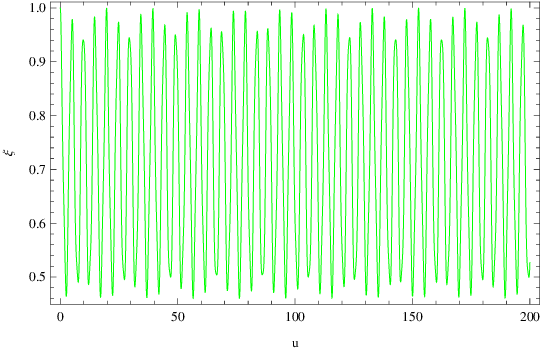}
			\caption{3D phase portrait}
			\label{fig6b}
		\end{subfigure}
		\begin{subfigure}{0.45\textwidth}
			\centering
			\includegraphics[width=.85\linewidth]{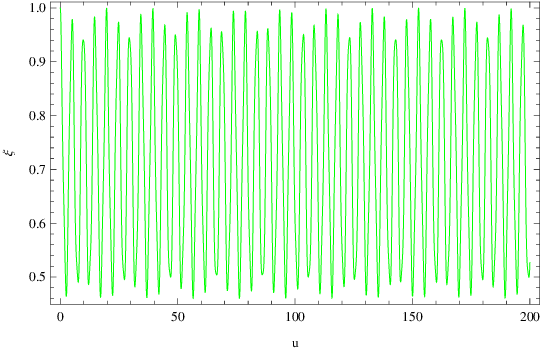}
			\caption{Time series }
			\label{fig6c}
		\end{subfigure}
		\caption{Phase portrait of system(\ref{3.10}) with different  $\sigma = -1.5$, $\alpha = 1$, $\omega = -3.5$, $z_0 =
			-0.3$, $\nu = -3.5$
			. }
		\label{fig6:main}
	\end{figure}

	Here, we visualized the dynamics of the system (\ref{3.10}) by varying the parameters. In fig.(\ref{fig4:main}) we display the 2D, 3D phase portrait and time series of the system \ref{3.10} by taking $\sigma = -1.1$, $\alpha = 1$, $\omega = 3.9$, $z_0 =
	-1.3$, $\nu = -4.5$. Here it is clearly visible that the system (\ref{3.10}) is a quasi-periodic system. Fig.(\ref{fig5:main}) demonstrate the 2D, 3D phase portrait and time series of the system (\ref{3.10}) with  $\sigma = -1.05$, $\alpha = 1$, $\omega = -3.5$, $z_0 =
	-0.3$, $\nu = -2.5$ as a result that the system (\ref{3.10}) is a periodic system by small changes of parameter.  Fig.(\ref{fig6:main}) display the the 2D, 3D phase portrait and time series of the system (\ref{3.10}) with  $\sigma = -1.5$, $\alpha = 1$, $\omega = -3.5$, $z_0 =
	-0.3$, $\nu = -3.5$ as a result that the system (\ref{3.10}) is a quasi-periodic system by a small perturbed term and frequency.
	

	\subsection{Sensitive analysis:}
	The sensitive analysis illustrates the system's dynamics that are highly dependent on the variation of parameters and initial conditions. In this section, we see how a system is highly sensitive to its initial condition. It is also a dynamic approach to ascertain whether the system is chaotic or not, because an essential condition of chaos is sensitive to initial condition.\\
	
\begin{figure}[htbp]
		\centering
		
		\begin{subfigure}{0.45\textwidth}
			\centering
			\includegraphics[width=\linewidth]{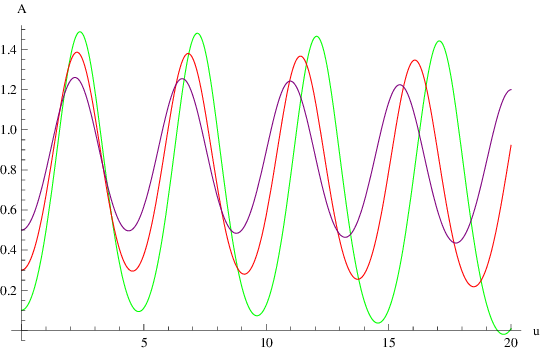}
			\caption{Green curve exhibit the initial condition $(u,A)=(0.1,0 )$, red curve exhibit the initial condition $(u,A)=(0.3,0 )$ and purple curve exhibit the initial condition $(u,A)=(0.5,0 )$ }
			\label{fig7a}
		\end{subfigure}
		\hfill
		\begin{subfigure}{0.45\textwidth}
			\centering
			\includegraphics[width=\linewidth]{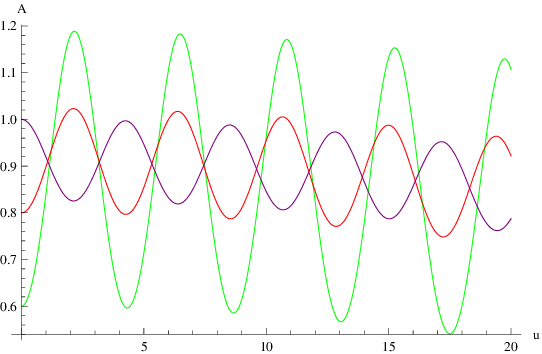}
			\caption{Green curve exhibit the initial condition $(u,A)=(0.6,0 )$, red curve exhibit the initial condition $(u,A)=(0.8,0 )$ and purple curve exhibit the initial condition $(u,A)=(1,0 )$}
			\label{fig7b}
		\end{subfigure}
		\caption{Sensitive analysis of (\ref{3.10}) with $\sigma = -2.5$, $\alpha = 5.5$, $\omega =- 3.5$, $z_0 =
			0.3$, $\nu = -0.05$. }
		\label{fig7:main}
	\end{figure}

	Fig (\ref{fig7a}) and Fig(\ref{fig7b}) illustrate the sensitivity to initial condition, which reveals the indication of chaos of the system (\ref{3.10}). The  green curve exhibits the initial condition $(u,A)=(0.1,0 )$, the red curve exhibit the initial condition $(u,A)=(0.3,0 )$ and the purple curve exhibit the initial condition $(u,A)=(0.5,0 )$ of Fig (\ref{fig7a}). The green curve exhibits the initial condition $(u,A)=(0.6,0 )$, red curve exhibit the initial condition $(u,A)=(0.8,0 )$ and the purple curve exhibits the initial condition $(u,A)=(1,0 )$ of Fig (\ref{fig7b}). Here, we can see that a simple modification of the initial condition gives a different solution of the system (\ref{3.10}).
	
	\subsection{Chaotic attractor:}
\begin{figure}[htbp]
		\centering
		
		\begin{subfigure}{0.45\textwidth}
			\centering
			\includegraphics[width=\linewidth]{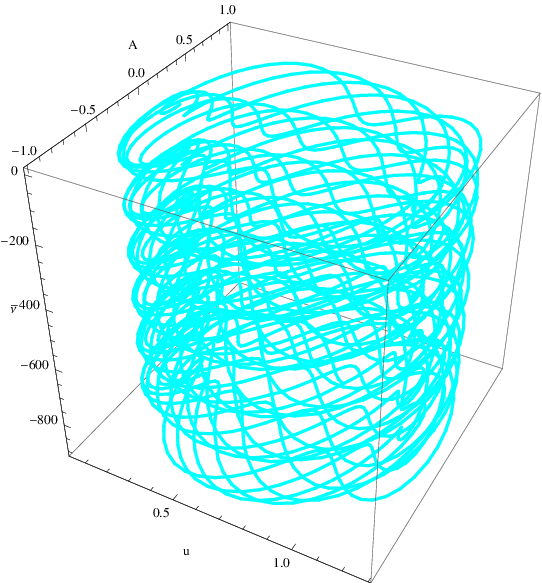}
			\caption{ }
			\label{fig8a}
		\end{subfigure}
		\hfill
		\begin{subfigure}{0.45\textwidth}
			\centering
			\includegraphics[width=\linewidth]{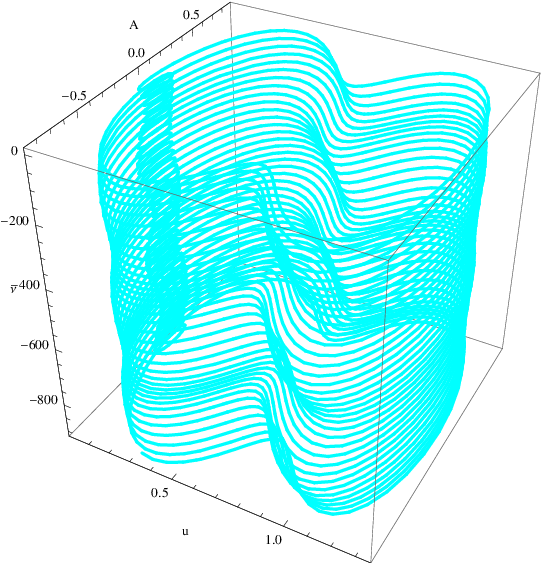}
			\caption{}
			\label{fig8b}
		\end{subfigure}
		\caption{Chaotic attractor of system(\ref{3.11}) with $\sigma = -2.5$, $\alpha = 1$, $\omega = 3.5$, $z_0 =
			-1.3$, $\nu = -4.5$.  }
		\label{fig8:main}
	\end{figure}

	A chaotic attractor can greatly aid comprehension of the system's long-term behavior in the chaotic system. It's fractal geometry and sensitive to initial conditions, revealing the presence of chaos. Fig.(\ref{fig8:main}) demonstrate the chaotic attractor, where fig.(\ref{fig8a}) portray the chaotic atrractor with initial condition $u(x,y,z)=(0.2,0,0.2)$ and fig.(\ref{fig8b}) portray the chaotic attractor with initial condition  $u(x,y,z)=(0.1,0,0.1)$ of the system (\ref{3.11}), that is highly sensitive to initial condition. The system (\ref{3.11}) is the improved version of the system (\ref{3.10}) which is given below,\\
	\begin{equation}\label{3.11}
		\begin{cases}
			\dfrac{du}{d\xi} = A, \\[6pt]
			\dfrac{dA}{d\xi} = p_{2}u^{2} + p_{1}u + p_{0}+Z_0 \cos(\nu \xi),  \\[6pt]
			\dfrac{d\nu}{d\xi}=C.
			\\[6pt]
		\end{cases}
	\end{equation}
	
	\subsection{Recurrence plot:}
\begin{figure}[h!]
		\centering
		\includegraphics[scale=0.4]{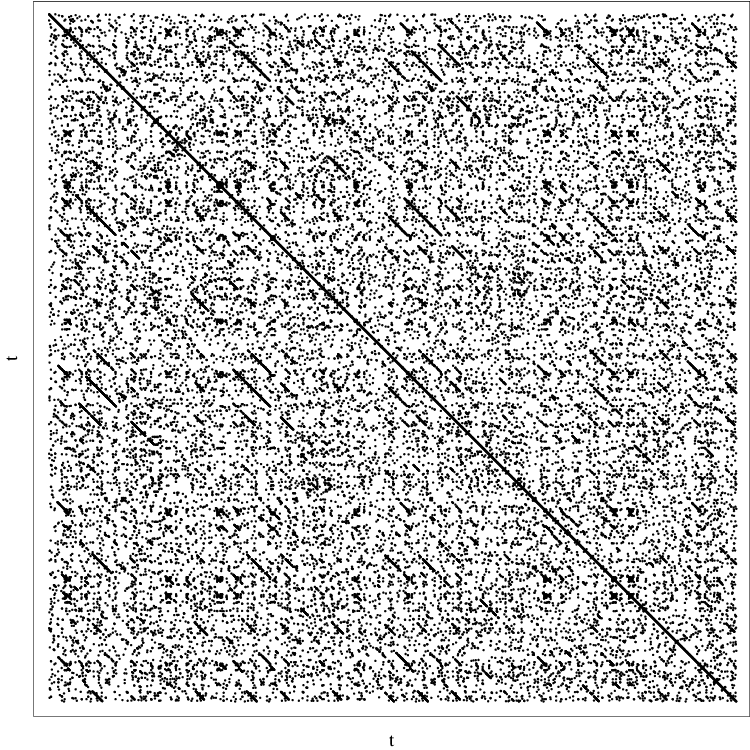}
		\caption{Recurrence of system(\ref{3.10}) with $p_0 = 0.7$, $p_1 = 0.2$, $z_0 =
			0.3$, $\nu = -4.5$.  }
		\label{fig9}
	\end{figure}

	A recurrence plot is a visual representation of the hidden dynamical structure of a nonlinear dynamical system to identify the behavior of the system within a time series. We can also say that it demonstrates the picture of a matrix where a black dot represents time $i$ and $j$ when a system returns to its previous state and a white space represents time  $i$ and $j$ when a system does not return to its previous state.\\
	Fig.(\ref{fig9}) illustrates the recurrence plot of the system (\ref{3.10}) with $p_0=0.7$, $p_1=0.2$, $Z_0=0.3$ and $\nu=-4.5$. The fragmented diagonal lines indicate the presence of chaos\cite{Eckmann1987}.
	
	\subsection{fractal dimension:}
\begin{figure}[h!]
		\centering
		\includegraphics[scale=0.5]{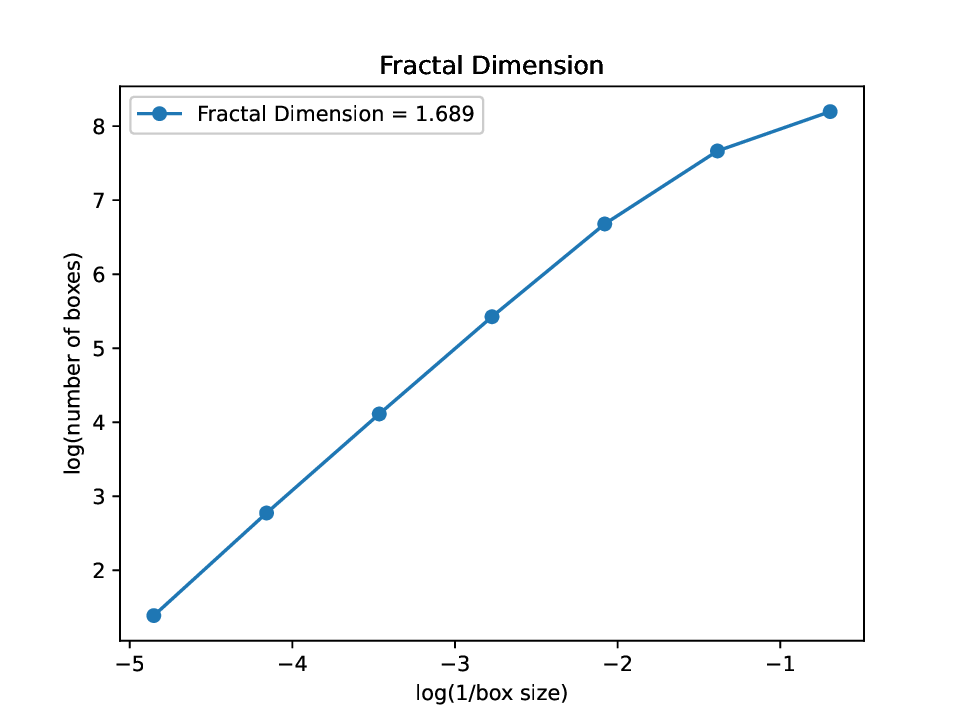}
		\caption{fractal dimension of the system(\ref{3.10}) with$p_0 = 0.7$, $p_1 = 0.2$, $z_0 =
			0.3$, $\nu = -4.5$.  }
		\label{fig10}
	\end{figure}

	The fractal dimension is a primary numerical indicator used to delineate the geometric complexity of an attractor that originates from a dynamical system. An integer values corresponding to the dimension of the limit cycle indicate the presence of periodicity, while a non-integer value provides insight into the complex geometric structure of phase space, which indicates the presence of chaos. The fractal dimension of the system (\ref{3.10}) is 1.689, which is a fractional number that indicates the presence of chaos. Fig.(\ref{fig10}) display the fractral dimension of system (\ref{3.10})\cite{Beenish2025}.

	\subsection{Power spectrum:}
\begin{figure}[htbp]
		\centering
		
		\begin{subfigure}{0.65\textwidth}
			\centering
			\includegraphics[width=\linewidth]{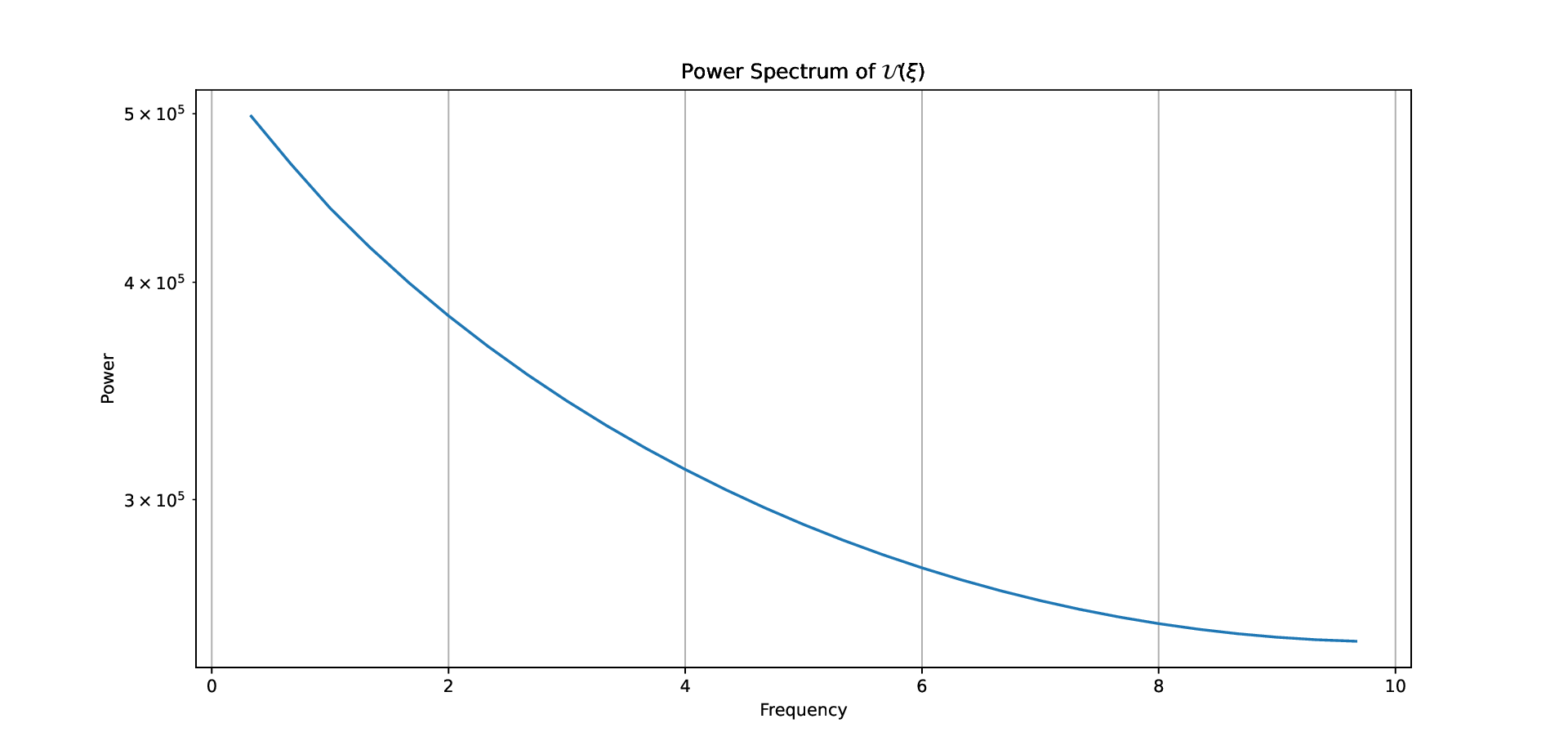}
			\caption{Power spectrum of the system(\ref{3.10}) with$p_0 = -0.7$, $p_1 = -0.5$, $z_0 =
				-3.3$, $\nu = -4.5$  that indicate the presence of chaos.  }
			\label{fig11a}
		\end{subfigure}
		\hfill
		\begin{subfigure}{0.65\textwidth}
			\centering
			\includegraphics[width=\linewidth]{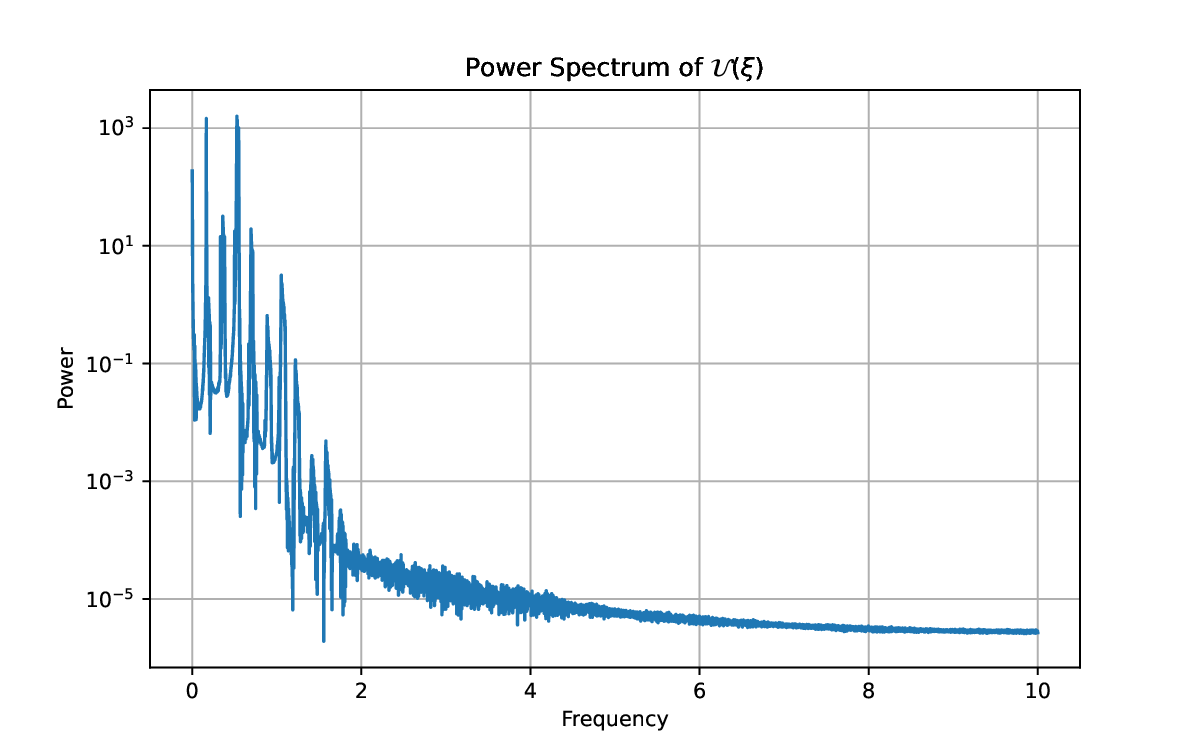}
			\caption{Power spectrum of the system(\ref{3.10}) with$p_0 = 18.5$, $p_1 = -6.5$, $z_0 =
				-5.3$, $\nu = -1.05$. that indicate the presence of of periodicity,}
			\label{fig11b}
		\end{subfigure}
		\caption{Power spectrum of the system(\ref{3.10}). }
		\label{fig11:main}
	\end{figure}

	The power spectrum is a qualitative tool that represents the frequency domain structure of a dynamical system. The power spectrum allows for the identification of discrepancies between periodic, quasi-periodic, and chaotic behavior of the system. A broadband continuous spectrum reflecting the presence of chaos, in contrast to a set of discrete but densely packed peaks at fundamental frequencies that indicate the presence of periodicity. Fig.(\ref{fig11:main}) represents the power spectrum of system(\ref{3.10})\cite{Almheidat2025}.
	\section{Conclusion:-}
	This study successfully investigates the  (2 + 1)-dimensional seventh-order Caudrey-Dodd-Gibbon-KP (sCDG-KP) equation, in which we employed the  ($\frac{G^{\prime}}{G^{\prime}+G+A}$) method to examine the exact solution of the Caudrey-Dod-Gibbon-KP (sCDG-KP) equation. Consequently, we visualized the 2D, 3D, and heat map of the obtained solutions. Next, our main focus goes to analyze the dynamical behavior of the system, and that system originates by altering the Caudrey-Dodd-Gibbon-KP (sCDG-KP) equation into a 2D system of equations. We use bifurcation analysis to determine the stability of the equilibrium point and to analyze how a parameter change causes changes in the system's long-term behavior. In this section, we also use different tools like chaotic attractor, sensitive analysis, fractal dimension, recurrence plot, and power spectrum to examine the chaotic behavior of the system; consequently, this allows us to reveal the complex behavior of the system and demonstrate how its solutions act under varying initial conditions and parameters. In future work, we will focus on expanding this method to fractional order PDE and use different tools like Lyapunov exponent, return map, and Poincar\'e map to identify the system's behavior.

	\section*{Conflicts of interest and funding}
	Conflict of interest: The authors state that they do not have any conflicts of interest.
	
	Funding: There is no financial assistance or funding available for the development of this manuscript.
	
	\section*{Data Availability Statement}
	On reasonable request, the corresponding author will make accessible the datasets created and/or analysed during the current work.
	\baselineskip18pt
	
	\bibliographystyle{plain}
	\bibliography{bib}
	
\end{document}